\documentclass[11pt]{article}
\usepackage{epic,latexsym,amssymb}
\usepackage{color}
\usepackage{tikz}
\usepackage{amsfonts,epsf,amsmath}

\newtheorem{thm}{Theorem}
\newtheorem{lem}{Lemma}

\newtheorem{cor}{Corollary}

\newtheorem{conj}{Conjecture}

\newcommand{\qed}{$\Box$}

\newcommand{\proof}{\noindent\textbf{Proof. }}

\let\oldenumerate\enumerate
\renewcommand{\enumerate}{
  \oldenumerate
  \setlength{\itemsep}{0pt}
  \setlength{\parskip}{0pt}
  \setlength{\parsep}{0pt}
}

\begin{document}

\title{Independence Number and Connectivity of Maximal Connected Domination Vertex Critical Graphs}

\author{Norah Almalki$^1$ \, and \, Pawaton Kaemawicharnurat$^{2,3}$
\\ \\
$^{1}$Department of Mathematics and Statistics \\
College of Science, Taif University, Saudi Arabia\\
\small \tt Email: norah@tu.edu.sa \\
\\
$^{2}$Department of Mathematics, Faculty of Science\\
King Mongkut's University of Technology Thonburi, Thailand\\
$^{3}$Mathematics and Statistics with Applications~(MaSA) \\
Bangkok, Thailand \\
\small \tt Email: pawaton.kae@kmutt.ac.th}

\maketitle

\begin{abstract}
A graph $G$ is said to be $k$-$\gamma_{c}$-edge critical if the connected domination number $\gamma_{c}(G) = k$ and $\gamma_{c}(G + uv) < k$ for every $uv \in E(\overline{G})$. A $2$-connected graph $G$ is said to be $k$-$\gamma_{c}$-vertex critical if $\gamma_{c}(G) = k$ and $\gamma_{c}(G - v) < k$ for any $v \in V(G)$. A maximal $k$-$\gamma_{c}$-vertex critical graph is a graph which are both $k$-$\gamma_{c}$-edge critical and $k$-$\gamma_{c}$-vertex critical. Let $\delta, \kappa$ and $\alpha$ be respectively the minimum degree, the connectivity and the independence number. In this paper, we prove that if $G$ is a maximal $3$-$\gamma_{c}$-vertex critical graph satisfying $\alpha = \kappa$, then the equation $\kappa = \delta$ still holds. We provide a class of maximal $3$-$\gamma_{c}$-vertex critical graphs with $\alpha < \kappa$ and $\kappa < \delta$. We finally show that every two vertices of any $3$-connected maximal $3$-$\gamma_{c}$-vertex critical graph are joined by hamiltonian path when $\kappa < \delta$.
\end{abstract}


\section{\bf Introduction}
Our basic graph theoretic notation and terminology follows that of Bondy and Murty\cite{BM}. Thus $G$ denotes a finite graph with vertex set $V(G)$ and edge set $E(G)$. For $S \subseteq V(G)$, $G[S]$ denotes the subgraph of $G$ \emph{induced} by $S$. Throughout this paper all graphs are simple and connected. The \emph{open neighborhood} $N_{G}(v)$ of a vertex $v$ in $G$ is $\{u \in V(G) : uv \in E(G)\}$. Further, the \emph{closed neighborhood} $N_{G}[v]$ of a vertex $v$ in $G$ is $N_{G}(v) \cup \{v\}$. The \emph{degree} $deg_{G}(v)$ of a vertex $v$ in $G$ is $|N_{G}(v)|$. The \emph{minimum degree} of a graph $G$ is denoted by $\delta(G)$. $N_{S}(v)$ denotes $N_{G}(v) \cap S$ where $S$ is a vertex subset of $G$. A \emph{tree} is a connected graph contains no cycle. A \emph{star} $K_{1, n}$ is a tree containing $n$ vertices of degree $1$ and exactly one vertex of degree $n$. An \emph{independent set} is a set of pairwise non-adjacent vertices. The \emph{independence number} $\alpha(G)$ is the maximum cardinality of an independent set. For a connected graph $G$, a \emph{cut set} is a vertex subset $S \subseteq V(G)$ such that $G - S$ is connected. The minimum cardinality of a vertex cut set of a graph $G$ is called the \emph{connectivity} and is denoted by $\kappa(G)$. If $G$ has $S = \{a\}$ as a minimum cut set, then $a$ is a \emph{cut vertex} of $G$ and $\kappa(G) = 1$. A graph $G$ is \emph{$s$-connected} if $\kappa(G) \geq s$. When ambiguity occur, we abbreviate $\delta(G), \alpha(G)$ and $\kappa(G)$ to $\delta, \alpha$ and $\kappa$, respectively. A \emph{hamiltonian path} is a path containing every vertex of a graph. A graph $G$ is \emph{hamiltonian connected} if every two vertices of $G$ are joined by a hamiltonian path. Obviously, a graph of connectivity one is not hamiltonian connected. Moreover, a graph of connectivity two does not contain a hamiltonian path between the two vertices in a minimum cut set. Thus we always focus on $3$-connected graphs when we study on the hamiltonian connected property of graphs. For a graph $G$, \emph{the Mycielskian} $\mu(G)$ of $G$ is the graph with vertex set $V(G) \cup V' \cup \{x\}$ where $V' = \{u'|u \in V(G)\}$ and edge set $E(G) \cup \{uv'|uv \in E(G)\} \cup \{v'x|v' \in V'\}$.
\vskip 5 pt

\indent For subsets $D, X \subseteq V(G)$, $D$ \emph{dominates} $X$ if every vertex in $X$ is either in $D$ or adjacent to a vertex in $D$. If $D$ dominates $X$, then we write $D \succ X$, further, we write $a \succ X$ when $D = \{a\}$. If $X = V(G)$, then $D$ is a \emph{dominating set} of $G$ and we write $D \succ G$ instead of $D \succ V(G)$. A \emph{connected dominating set} of a graph $G$ is a dominating set $D$ of $G$ such that $G[D]$ is connected. If $D$ is a connected dominating set of $G$, we then write $D \succ_{c} G$. A smallest connected dominating set is call a $\gamma_{c}$\emph{-set}. The cardinality of a $\gamma_{c}$-set of $G$ is called the \emph{connected domination number} of $G$ and is denoted by $\gamma_{c}(G)$. A \emph{total dominating set} of a graph $G$ is a subset $D$ of $V(G)$ such that every vertex in $G$ is adjacent to a vertex in $D$. The minimum cardinality of a total dominating set of $G$ is called the \emph{total domination number} of $G$ and is denoted by $\gamma_{t}(G)$.
\vskip 5 pt

\indent A graph $G$ is $k$-$\gamma_{c}$-\emph{edge critical} if $\gamma_{c}(G) = k$ and $\gamma_{c}(G + uv) < k$ for each pair of non-adjacent vertices $u$ and $v$ of $G$. A graph $G$ is $k$-$\gamma_{c}$-\emph{vertex critical} if $\gamma_{c}(G) = k$ and $\gamma_{c}(G - v) < k$ for any vertex $v$ of $G$. A graph $G$ is a \emph{maximal} $k$-$\gamma_{c}$-\emph{vertex critical} if $G$ is both $k$-$\gamma_{c}$-edge critical and $k$-$\gamma_{c}$-vertex critical. It is easy to see that a disconnected graph cannot contain a connected dominating set. Thus if a graph $G$ contains a cut vertex $c$, then $G - c$ has no connected dominating set. Hence, in the following, we always assume that $k$-$\gamma_{c}$-vertex critical graphs are $2$-connected. A $k$-$\gamma_{t}$-\emph{edge critical} graph is similarly defined.
\vskip 12 pt

\indent This paper focuses on the relationship of the connectivity and the independence number of $3$-$\gamma_{c}$-edge critical graphs. For related results on $3$-$\gamma$-edge critical graphs, Zhang and Tian\cite{ZT} showed that the independence number of every $3$-$\gamma$-edge critical graph does not exceed $\kappa + 2$ and, moreover if $\alpha = \kappa + 2$, then $\kappa = \delta$. They used this result to established, the short proof, that every $3$-$\gamma$-edge critical graph contains a hamiltonian path. Further, Kaemawichanurat and Caccetta \cite{KC} showed that every $3$-$\gamma_{c}$-edge critical graph satisfies $\alpha \leq \kappa + 2$. Moreover, if $G$ is a $3$-$\gamma_{c}$-edge critical graph such that $\kappa + 1 \leq \alpha \leq \kappa + 2$, then $\kappa = \delta$ with only one exception. They also showed that if $\kappa \geq 3$, then $\alpha = \kappa + p$ if and only if $\alpha = \delta + p$ for all $p \in \{1, 2\}$. The condition $\kappa + 1 \leq \alpha \leq \kappa + 2$ is best possible to establish that $\kappa = \delta$.
\vskip 5 pt

\indent In this paper, we prove that if a $3$-$\gamma_{c}$-edge critical graph $G$ is $3$-$\gamma_{c}$-vertex critical ($G$ is maximal $3$-$\gamma_{c}$-vertex critical) satisfying $\alpha = \kappa$, then $\kappa = \delta$. We provide a class of maximal $3$-$\gamma_{c}$-vertex critical graphs with $\alpha < \kappa$ and $\kappa < \delta$. Hence, the condition $\alpha = \kappa$ is needed to prove that $\kappa = \delta$ in maximal $3$-$\gamma_{c}$-vertex critical graphs. We finish this work by showing that every $3$-connected maximal $3$-$\gamma_{c}$-vertex critical graph are hamiltonian connected if $\kappa < \delta$.
\vskip 5 pt

\section{\bf Preliminaries}
\noindent In this section, we state a number of results that we make use of in establishing our theorems. We begin with a result of Chv\'{a}tal and Erd\"{o}s\cite{CE} which gives hamiltonian properties of graphs according to the connectivity and the independence number.
\vskip 5 pt

\begin{thm}\cite{CE}\label{thm CE}
Let $G$ be an $s$-connected graph containing no independent set of $s$ vertices. Then $G$ is hamiltonian connected.
\end{thm}
\vskip 5 pt

\indent On $3$-$\gamma_{c}$-edge critical graphs, Ananchuen\cite{A} characterized such graphs with a cut vertex. For positive integers $n_{i}$ and $r \geq 2$, let $H = \cup^{r}_{i = 1}K_{1, n_{i}}$. For $1 \leq j \leq r$, let $c_{j}$ be the center of $K_{1, n_{j}}$ in $H$ and $w^{j}_{1}, w^{j}_{2}, ..., w^{j}_{n_{j}}$ the end vertices of $K_{1, n_{j}}$ in $H$. The graphs $G_{1}$ and $G_{2}$ are defined as follows. Set $V(G_{1}) = V(H) \cup \{x, y\}$ and $E(G_{1}) = E(\overline{H}) \cup \{xy\} \cup \{xw^{j}_{i}: 1 \leq i \leq n_{j}$ and $1 \leq j \leq r\}$. Set $V(G_{2}) = V(H) \cup \{x, y\} \cup U$ where $|U| \geq 1$ and $E(G_{2}) = E(\overline{H}) \cup \{xy\} \cup \{xw^{j}_{i}: 1 \leq i \leq n_{j}$ and $1 \leq j \leq r\} \cup \{uz: u \in U$ and $z \in V(G_{1}) - \{x, y, u\}\}$.



\indent The following two lemmas, observed by Chen et al.\cite{CSM}, give fundamental properties of $3$-$\gamma_{c}$-edge critical graphs.
\vskip 5 pt

\begin{lem}\label{lem 1}\cite{CSM}
Let $G$ be a $3$-$\gamma_{c}$-edge critical graph and, for a pair of non-adjacent vertices $u$ and $v$ of $G$, let $D_{uv}$ be a $\gamma_{c}$-set of $G + uv$. Then
\begin {enumerate}
\item [(1)] $|D_{uv}| = 2$,
\item [(2)] $D_{uv} \cap \{u, v\} \neq \emptyset$ and
\item [(3)] if $u \in D_{uv}$ and $v \notin D_{uv}$, then $N_{G}(v) \cap D_{uv} = \emptyset$.
\end {enumerate}
\end{lem}
\vskip 5 pt

\begin{lem}\label{lem 2}\cite{CSM}
Let $G$ be a $3$-$\gamma_{c}$-edge critical graph and $I$ an independent set with $|I| = p \geq 3$. Then the vertices in $I$ can be ordered as $a_{1}, a_{2}, ..., a_{p}$ and there exists a path $x_{1}, x_{2},..., x_{p - 1}$ in $G - I$ with $\{a_{i}, x_{i}\} \succ_{c} G + a_{i}a_{i + 1}$ for $1 \leq i \leq p - 1$.
\end{lem}
\vskip 5 pt

\indent On $3$-$\gamma_{c}$-vertex critical graphs, Ananchuen et al.\cite{AAP} provided fundamental properties of such graphs.
\vskip 5 pt

\begin{lem}\label{lem 6}\cite{AAP}
Let $G$ be a $3$-$\gamma_{c}$-vertex critical graph and, for a vertex $v$ of $G$, let $D_{v}$ be a $\gamma_{c}$-set of $G - v$. Then
\begin {enumerate}
\item [(1)] $|D_{v}| = 2$ and
\item [(2)] $D_{v} \cap N_{G}[v] = \emptyset$.
\end {enumerate}
\end{lem}
\vskip 5 pt

\indent In the study of $3$-$\gamma_{t}$-edge critical graphs, Simmons\cite{S} proved that such graphs satisfy $\alpha \leq \delta + 2$. It was pointed out in Ananchuen\cite{A} that a $3$-$\gamma_{c}$-edge critical graph is also $3$-$\gamma_{t}$-edge critical and vice versa. Then every $3$-$\gamma_{c}$-edge critical graph satisfies $\alpha \leq \delta + 2$. In \cite{KCA}, we established the result of $3$-$\gamma_{c}$-edge critical graphs when $\alpha = \delta + 2$.
\vskip 5 pt

\begin{thm}\label{thm K}(\cite{S}, \cite{KCA})
Let $G$ be a $3$-$\gamma_{c}$-edge critical graph with $\delta \geq 2$. Then $\alpha \leq \delta + 2$. Moreover if  $\alpha = \delta + 2$, then $G$ contains exactly one vertex $x$ of degree $\delta$ and $G[N[x]]$ is a clique.
\end{thm}
\vskip 5 pt

\indent In \cite{KCA1}, we established results of maximal $3$-$\gamma_{c}$-vertex critical graphs. Let $G$ be a maximal $3$-$\gamma_{c}$-vertex critical graph with a vertex cut set $S$ and $C_{1}, C_{2}, ..., C_{m}$ be the components of $G - S$.
\vskip 5 pt

\begin{lem}\label{lem Co1p}\cite{KCA1}
For all $v \in V(G)$, if $m \geq 3$ or $v \in S \cup V(C_{i})$ where $|V(C_{i})| > 1$, then $G$ satisfies these following properties.
\begin {enumerate}
\item [(1)] $D_{v} \cap S \neq \emptyset$.
\item [(2)] $v$ does not dominate $S$.
\end {enumerate}
\end{lem}
\vskip 5 pt

\begin{lem}\label{lem t1p}\cite{KCA1}
Let $a \in V(C_{i})$ for some $i \in \{1, 2, ..., m\}$. Then $G$ has these following properties.
\begin {enumerate}
\item [(1)] Let $b \in V(C_{j})$ for some $j \in \{1, 2, ..., m\}$ such that $\{a, b\}$ does not dominate $G$, if $m \geq 3$ or $|V(C_{i})|, |V(C_{j})| > 1$, then $|D_{ab} \cap \{a, b\}| = 1$ and $|D_{ab} \cap S| = 1$.
\item [(2)] If $c \in D_{a}$ where $c$ is an isolated vertex in $S$, then $m = 2$ and $\{w\} = V(C_{j})$ for some $j \in \{1, 2\}$ where $\{w\} = D_{a} - \{c\}$.
\end {enumerate}
\end{lem}
\vskip 5 pt

\indent We characterized all maximal $3$-$\gamma_{c}$-vertex critical graphs of connectivity $s \geq 2$ with a smallest cut set contains no edge
\vskip 5 pt

\begin{thm}\label{thm C3p}\cite{KCA1}
Let $G$ be a maximal $3$-$\gamma_{c}$-vertex critical graph of connectivity $s \geq 2$ with a smallest cut set contains no edge. Then $G$ is isomorphic to $G_{3} = \mu(K_{s})$.
\end{thm}
\vskip 5 pt

\setlength{\unitlength}{0.7cm}
\begin{center}
\begin{picture}(10,5)

\put(1, 3){\circle*{0.2}}
\put(3, 3){\circle*{0.2}}
\put(3, 4){\circle*{0.2}}
\put(3, 5){\circle*{0.2}}
\put(3, 1){\circle*{0.2}}

\put(6, 3){\circle*{0.2}}
\put(6, 4){\circle*{0.2}}
\put(6, 5){\circle*{0.2}}
\put(6, 1){\circle*{0.2}}

\put(3, 1.7){\circle*{0.08}}
\put(3, 2){\circle*{0.08}}
\put(3, 2.3){\circle*{0.08}}

\multiput(3, 3)(0.4,0){8}{\line(1,0){0.2}}
\multiput(3, 4)(0.4,0){8}{\line(1,0){0.2}}
\multiput(3, 5)(0.4,0){8}{\line(1,0){0.2}}
\multiput(3, 1)(0.4,0){8}{\line(1,0){0.2}}

\put(6, 3){\oval(1, 5)}

\put(1, 3){\line(1, 1){2}}
\put(1, 3){\line(2, 1){2}}
\put(1, 3){\line(1, 0){2}}
\put(1, 3){\line(1, -1){2}}

\put(3, 5){\line(3, -1){3}}
\put(3, 5){\line(3, -2){3}}
\put(3, 5){\line(3, -4){3}}

\put(3, 4){\line(3, 1){3}}
\put(3, 4){\line(1, -1){3}}
\put(3, 4){\line(3, -1){3}}

\put(3, 3){\line(3, -2){3}}
\put(3, 3){\line(3, 2){3}}
\put(3, 3){\line(3, 1){3}}

\put(3, 1){\line(3, 2){3}}
\put(3, 1){\line(3, 4){3}}
\put(3, 1){\line(1, 1){3}}

\put(0.8, 3.3){\footnotesize$v$}

\put(1, -0.5){\footnotesize\textbf{Figure 1 :} A graph $G_{3} = \mu(K_{s})$}

\end{picture}
\end{center}
\vskip 30 pt

\indent In \cite{KCA2}, we established the upper bound of the independence number of maximal $3$-$\gamma_{c}$-vertex critical graphs in term of the minimum degree.
\begin{thm}\label{thm I3p}\cite{KCA2}
Let $G$ be a maximal $3$-$\gamma_{c}$-vertex critical graph. Then $\alpha \leq \delta$.
\end{thm}

\section{Maximal $3$-$\gamma_{c}$-Vertex Critical Graphs}
In this section we show, by using Theorem \ref{thm I3p}, that every maximal $3$-$\gamma_{c}$-vertex critical graph satisfies $\alpha \leq \kappa$. We also study the such graphs when $\alpha = \kappa$. Since we characterized all maximal $3$-$\gamma_{c}$-vertex critical graphs with $\kappa \leq 3$ in \cite{KCA1}, in the following, we focus on $|S| = \kappa \geq 4$. We define $\alpha_{1}, \alpha_{2}, p, S, H_{1}$ and $H_{2}$ by the same as in the previous section.
\vskip 5 pt

\begin{thm}\label{thm Co2p}
The independence number of any maximal $3$-$\gamma_{c}$-vertex critical graph does not exceed the connectivity.
\end{thm}
\proof Suppose to the contrary that $\kappa + 1 \leq \alpha$. So $|S| + 1 \leq \alpha_{1} + \alpha_{2} + |S \cap I|$. Then
\begin{align}\label{eq 1}
|S - I| + 1  & = |S| - |S \cap I| + 1 \leq \alpha_{1} + \alpha_{2}
\end{align}

\noindent \textbf{Claim 1 :} $|V(C_{i})| > 1$ for all $i \in \{1, 2, ..., m\}$, in particular, $|H_{i}| > 1$.\\
\indent Suppose $V(C_{i}) = \{c\}$ for some $i \in \{1, 2, ..., m\}$. Thus $N_{G}(c) \subseteq S$. By Theorem \ref{thm I3p},
\begin{center}
  $\delta \leq deg_{G}(c) \leq |S| < |S| + 1 = \kappa + 1 \leq \alpha \leq \delta$,
\end{center}
\noindent a contradiction and thus establishing Claim 1.
\vskip 5 pt

\indent Let $p = \alpha_{1} + \alpha_{2}$ and $\{a_{1}, a_{2}, ..., a_{p}\} = \cup^{2}_{i = 1}I_{i}$. If $p = 1$, then, by (\ref{eq 1}), $|S - I| = 0$. It follows that $S = S \cap I$. Thus $S$ is an independent set. By Theorem \ref{thm C3p}, $G$ is $G_{3}$. Hence, $N_{G_{3}}(x)$ in the graph $G_{3}$ is a minimum cut set such that $G_{3} - N_{G_{3}}(x)$ contains $x$ as a singleton component, contradicting Claim 1. Therefore $p \geq 2$.
\vskip 5 pt

\noindent \textbf{Claim 2 :} For $a, b \in \cup^{2}_{i = 1}I_{i}$, $|D_{ab} \cap \{a, b\}| = 1$ and $|D_{ab} \cap (S - I)| = 1$.\\
\indent Since $|S| \geq 4$ and $2 \leq p = \alpha_{1} + \alpha_{2}$. If $p \geq 3$, then $\cup^{2}_{i = 1}I_{i} - \{a, b\} \neq \emptyset$. If $p = 2$, then, by Equation \ref{eq 1}, $|S| - |S \cap I| + 1 \leq 2$. As $|S| \geq 4$, we must have $|S \cap I| \geq 3$, in particular, $S \cap I \neq \emptyset$. Hence $(\cup^{2}_{i = 1}I_{i} - \{a, b\}) \cup (S \cap I) \neq \emptyset$. So $\{a, b\}$ does not dominate $G$. By Claim 1 and Lemma \ref{lem t1p}(1), $|D_{ab} \cap \{a, b\}| = 1$ and $|D_{ab} \cap S| = 1$. Without loss of generality, let $a \in D_{ab}$ and $\{a'\} = D_{ab} \cap S$. By the connectedness of $(G + ab)[D_{ab}]$, $a' \in S - I$ and thus establishing Claim 2.
\vskip 5 pt

\indent Suppose that $p = 2$. Consider $G + a_{1}a_{2}$. By Claim 2, $|D_{a_{1}a_{2}} \cap (S - I)| = 1$. Since $D_{a_{1}a_{2}} \cap (S - I) \subseteq S - I$, by Equation \ref{eq 1},
\begin{center}
$1 \leq |S - I| \leq \alpha_{1} + \alpha_{2} - 1 = p - 1 = 1$.
\end{center}
\noindent Therefore, $D_{a_{1}a_{2}} \cap (S - I) = S - I$. If $p \geq 3$, then, by Lemma \ref{lem 2}, the vertices $a_{1}, a_{2}, ..., a_{p}$ can be ordered as $x_{1}, x_{2}, ..., x_{p}$ and there exists a path $y_{1}, y_{2}, ..., y_{p - 1}$ such that $\{x_{i}, y_{i}\} \succ_{c} G + x_{i}x_{i + 1}$ for $i \in \{1, 2, ..., p - 1\}$. Since $\{x_{1}, x_{2}, ..., x_{p}\} \subseteq \cup^{2}_{i = 1}I_{i}$, it follows by Claim 2 that $\{y_{1}, y_{2}, ..., y_{p - 1}\} \subseteq S - I$. Hence, by (\ref{eq 1}), $p - 1 \leq |S - I| \leq \alpha_{1} + \alpha_{2} - 1 = p - 1$. In both cases $p = 2$ and $p \geq 3$, we have that $\{y_{1}, y_{2}, ..., y_{p - 1}\} = S - I$.

\indent If $p = 2$, then clearly $G[\{y_{1}\}]$ is a clique. Suppose $p \geq 3$. Consider $G + x_{i}x_{j}$ for $2 \leq i \neq j \leq p$. By Claim 2, $|D_{x_{i}x_{j}} \cap \{x_{i}, x_{j}\}| = 1$ and $|D_{x_{i}x_{j}} \cap (S - I)| = 1$. Without loss of generality, let $x_{i} \in D_{x_{i}x_{j}}$. Because $S - I = \{y_{1}, y_{2}, ..., y_{p - 1}\}$, by Lemma \ref{lem 1}(3), $D_{x_{i}x_{j}} \cap (S - I) = \{y_{j - 1}\}$. Since $x_{i}y_{i - 1} \notin E(G)$, $y_{i - 1}y_{j - 1} \in E(G)$. Therefore, $G[\{y_{1}, y_{2}, ..., y_{p - 1}\}]$ is a clique. Since $\{x_{1}, x_{2}, ..., x_{p}\} \subseteq I$, $y_{i} \succ (S \cap I)$ for $i \in \{1, 2, ..., p - 1\}$. Therefore $y_{i} \succ S$, contradicting Lemma \ref{lem Co1p}(2). Hence, $\alpha \leq \kappa$.
\qed

\indent We have by Theorem \ref{thm C3p} that the graph $G_{3} = \mu(K_{s})$ has $N_{G_{3}}(x)$ as a minimum cut set and also a maximum independent set. Hence $\alpha(G_{3}) = \kappa(G_{3})$. Therefor, the bound in Theorem \ref{thm Co2p} is best possible. We, further, focus on maximal $3$-$\gamma_{c}$-vertex critical graphs which satisfy $\alpha = \kappa$. We then have $|S - I| + |S \cap I| = |S| = \alpha_{1} + \alpha_{2} + |S \cap I|$. So
\begin{align}\label{eq 2}
|S - I| = \alpha_{1} + \alpha_{2} = p.
\end{align}

\indent We may assume without loss of generality that $\alpha_{1} \leq \alpha_{2}$. We would like to show that if a maximal $3$-$\gamma_{c}$-vertex critical graph $G$ satisfies $\alpha = \kappa$, then $G - S$ contains at least one singleton component for every minimum cut set $S$. We then suppose to the contrary that there is no singleton component of $G - S$, in particular, $|H_{i}| > 1$ for $i = 1, 2$. We next establish the following lemmas.
\vskip 5 pt

\begin{lem}\label{lem Co3p}
Let $G$ be a maximal $3$-$\gamma_{c}$-vertex critical graph. If $\alpha = \kappa$ and $|V(C_{i})| > 1$ for all $i \in \{1, 2, ..., m\}$, then $p \geq 3$.
\end{lem}
\proof
By the assumption, $|H_{i}| > 1$ for $i = 1, 2$. Suppose first that $p = 0$. Thus $S = S \cap I$. By Theorem \ref{thm C3p}, $G$ is $G_{3}$. We have $N_{G_{3}}(x)$ is a minimum cut set of $G_{3}$ and $x$ is a singleton component of $G - N_{G_{3}}(x)$, contradicting the assumption. We distinguish $2$ cases.
\vskip 5 pt

\noindent \textbf{Case 1 :} $p = 1$.\\
By  Equation \ref{eq 2}, $|S - I| = 1$. Let $\{v\} = S - I, \{a_{1}\} = \cup^{2}_{i = 1}I_{i}$ and $\{a_{2}, a_{3}, ..., a_{\alpha}\} = S \cap I$. Therefore $\alpha_{1} = 0, \alpha_{2} = 1$. Therefore, $a_{1} \in H_{2}$. Since $|S| \geq 4$, $|S \cap I| \geq 3$. Lemma \ref{lem 2} yields that the vertices in $\{a_{2}, a_{3}, ..., a_{\alpha}\}$ can be ordered as $x_{1}, x_{2}, ..., x_{\alpha - 1}$ and there exists a path $P = y_{1}, y_{2}, ..., y_{\alpha - 2}$ with $\{x_{i}, y_{i}\} \succ_{c} G + x_{i}x_{i + 1}$ for $i \in \{1, .., \alpha - 2\}$. Since each $y_{i}$ is adjacent to at least one vertex of $I$ for $i = 1, 2, ...,  \alpha - 2$, $y_{i} \neq a_{1}$. To dominate $a_{1}$, $y_{i} \in H_{2} \cup \{v\}$.
\vskip 5 pt

\noindent \textbf{Subcase 1.1 :} $v \notin V(P)$.\\
So $V(P) \subseteq H_{2}$. Thus $x_{i} \succ H_{1}$ for $i = 1, 2, ...,  \alpha - 2$. Since $S$ is a minimum cut set, $N_{H_{1}}(v) \neq \emptyset$. Let $u \in N_{H_{1}}(v)$. So $u \succ \{x_{1}, x_{2}, ..., x_{\alpha - 2}, v\}$. Lemma \ref{lem Co1p}(2) implies that $ux_{\alpha - 1} \notin E(G)$. Consider $G + uy_{\alpha - 2}$. Since $ux_{\alpha - 1}, y_{\alpha - 2}x_{\alpha - 1} \notin E(G)$, by Lemma \ref{lem t1p}(1), $|D_{uy_{\alpha - 2}} \cap \{u, y_{\alpha - 2}\}| = 1$ and $|D_{uy_{\alpha - 2}} \cap S| = 1$. So either $y_{\alpha - 2} \in D_{uy_{\alpha - 2}}$ or $u \in D_{uy_{\alpha - 2}}$. In the first case, by Lemma \ref{lem 1}(3), $\{x_{1}, x_{2}, ..., x_{\alpha - 2}, v\} \cap D_{uy_{\alpha - 2}} = \emptyset$. Hence $x_{\alpha - 1} \in D_{uy_{\alpha - 2}}$. But $G[D_{uy_{\alpha - 2}}]$ is not connected. Thus $u \in D_{uy_{\alpha - 2}}$. By the connectedness of $(G + uy_{\alpha - 2})[D_{uy_{\alpha - 2}}]$, $x_{\alpha - 1} \notin D_{uy_{\alpha - 2}}$. If $x_{i} \in D_{uy_{\alpha - 2}}$ for $i \in \{1, 2, ..., \alpha - 2\}$, then $D_{uy_{\alpha - 2}}$ does not dominate $x_{\alpha - 1}$. Therefore $v \in D_{uy_{\alpha - 2}}$. Thus $va_{1} \in E(G)$. Consider $G + ua_{1}$. Since $ux_{\alpha - 1}, a_{1}x_{\alpha - 1} \notin E(G)$, by Lemma \ref{lem t1p}(1), $|D_{ua_{1}} \cap \{u, a_{1}\}| = 1$ and $|D_{ua_{1}} \cap S| = 1$. Hence either $u \in D_{ua_{1}}$ or $a_{1} \in D_{ua_{1}}$. In the first case, $v \notin D_{ua_{1}}$ by Lemma \ref{lem 1}(3). By the connectedness of $(G + ua_{1})[D_{ua_{1}}]$, $x_{\alpha - 1} \notin D_{ua_{1}}$. To dominate $x_{\alpha - 1}$, $D_{ua_{1}} \cap \{x_{1}, x_{2}, ..., x_{\alpha - 2}\} \neq \emptyset$. So $D_{ua_{1}} \cap S = \emptyset$, a contradiction. Hence, $a_{1} \in D_{ua_{1}}$. Lemma \ref{lem 1}(3) gives that $v \notin D_{ua_{1}}$. By the connectedness of $(G + ua_{1})[D_{ua_{1}}]$, $\{x_{1}, x_{2}, ..., x_{\alpha - 1}\} \cap D_{ua_{1}} = \emptyset$. Clearly $D_{ua_{1}} \cap S = \emptyset$, a contradiction and Case 1.1 cannot occur.
\vskip 5 pt

\noindent \textbf{Subcase 1.2 :} $v \in V(P)$.\\
Thus there exists $j \in \{1, 2, ..., \alpha - 2\}$ such that $y_{j} = v$. Thus $x_{i} \succ H_{1}$ for $i \neq j, \alpha - 1$ and $va_{1} \in E(G)$. Since $a_{1}, x_{\alpha - 1} \in I$, $a_{1}x_{\alpha - 1} \notin E(G)$. If $x_{\alpha - 1}$ is not adjacent to $w \in H_{1}$, then consider $G + wa_{1}$. Lemma \ref{lem t1p}(1) implies that $|D_{wa_{1}} \cap \{w, a_{1}\}| = 1$ and $|D_{wa_{1}} \cap S| = 1$. Thus either $w \in D_{wa_{1}}$ or $a_{1} \in D_{wa_{1}}$. In both case, $x_{\alpha - 1} \notin D_{wa_{1}}$ because $(G + wa_{1})[D_{wa_{1}}]$ is connected. If $w \in D_{wa_{1}}$, then, by Lemma \ref{lem 1}(3), $v \notin D_{wa_{1}}$. To dominate $x_{\alpha - 1}$, $\{x_{1}, x_{2}..., x_{\alpha - 2}\} \cap D_{wa_{1}} = \emptyset$. So $D_{wa_{1}} \cap S = \emptyset$, a contradiction. Hence, $a_{1} \in D_{wa_{1}}$.By the connectedness of $(G + wa_{1})[D_{wa_{1}}]$, $D_{wa_{1}} \cap \{x_{1}, x_{2}, ..., x_{\alpha - 1}\} = \emptyset$. To dominate $x_{j + 1}$, $v \notin D_{wa_{1}}$. We then have $D_{wa_{1}} \cap S = \emptyset$, a contradiction. Thus $x_{\alpha - 1} \succ H_{1}$. Clearly $x_{i} \succ H_{1}$ for $i \neq j$. Since $S$ is a minimum cut set, $N_{H_{1}}(v) \neq \emptyset$. Let $u' \in N_{H_{1}}(v)$. Lemma \ref{lem Co1p}(2) yields that $u' \succ S - \{x_{j}\}$. Consider $G + u'a_{1}$. By the same arguments as considering $G + ua_{1}$, we have a contradiction. Thus Case 1 cannot occur.
\vskip 5 pt

\noindent \textbf{Case 2 :} $p = 2$.\\
Let $\{a_{1}, a_{2}\} = \cup^{2}_{i = 1}I_{i}$. By Equation \ref{eq 2}, $|S - I| = p = 2$. Since $|S| \geq 4$, it follows that $|S \cap I| \geq 2$, in particular, $S \cap I \neq \emptyset$ and $\{a_{1}, a_{2}\}$ does not dominate $G$. Consider $G + a_{1}a_{2}$. Lemma \ref{lem t1p}(1) gives that $|D_{a_{1}a_{2}} \cap \{a_{1}, a_{2}\}| = 1$ and $|D_{a_{1}a_{2}} \cap S| = 1$. Without loss of generality, let $a_{1} \in D_{a_{1}a_{2}}$. By the connectedness of $(G + a_{1}a_{2})[D_{a_{1}a_{2}}]$, $|(S - I) \cap D_{a_{1}a_{2}}| = 1$. Let $\{u\} = (S - I) \cap D_{a_{1}a_{2}}$. Thus $ua_{1} \in E(G), ua_{2} \notin E(G)$ and $u \succ S \cap I$. Let $v \in S - (I \cup \{u\})$. Lemma \ref{lem Co1p}(2) implies that $uv \notin E(G)$. Therefore $a_{1}v \in E(G)$
\vskip 5 pt

\noindent \textbf{Subcase 2.1 :} $\alpha_{1} = 1$ and $\alpha_{2} = 1$.\\
Without loss of generality, let $a_{1} \in I_{1}$ and $a_{2} \in I_{2}$. Since $|S \cap I| \geq 2$, there exist $a_{3}, a_{4} \in S \cap I$. Consider $G + a_{3}a_{4}$. Lemma \ref{lem 1}(2) gives that $D_{a_{3}a_{4}} \cap \{a_{3}, a_{4}\} \neq \emptyset$. To dominate $a_{1}$, $D_{a_{3}a_{4}} \neq \{a_{3}, a_{4}\}$. Without loss of generality, let $a_{3} \in D_{a_{3}a_{4}}$. Lemma \ref{lem 1}(1) implies that $|D_{a_{3}a_{4}} - \{a_{3}\}| = 1$. Let $y \in D_{a_{3}a_{4}} - \{a_{3}\}$. To dominate $\{a_{1}, a_{2}\}$, $y \notin \cup^{2}_{i = 1}H_{i}$. By the connectedness of $(G + a_{3}a_{4})[D_{a_{3}a_{4}}]$, $y \in \{v, u\}$. Since $uv \notin E(G)$, $a_{3}u, a_{3}v \in E(G)$. Consider $G - a_{3}$ . Lemma \ref{lem 6}(2) then implies that $D_{a_{3}} \cap \{u, v\} = \emptyset$. Lemma \ref{lem Co1p}(1) implies also that $D_{a_{3}} \cap S \neq \emptyset$. So there exists $z \in D_{a_{3}} \cap (S \cap I)$. Lemma \ref{lem 6}(1) yields that $|D_{a_{3}} - \{z\}| = 1$. Let $D_{a_{3}} - \{z\} = \{z'\}$. Since $z \in S \cap I$, $za_{1} \notin E(G)$. Thus $z' \in H_{1}$ to dominate $a_{1}$. So $D_{a_{3}}$ does not dominate $a_{2}$, a contradiction and Subcase 2.1 cannot occur.
\vskip 5 pt

\noindent \textbf{Subcase 2.2 :} $\alpha_{1} = 0$ and $\alpha_{2} = 2$.\\
Thus $u \succ H_{1}$. Let $b_{1} \in H_{1}$. Clearly $\{a_{1}, b_{1}\}$ does not dominate $G$. Consider $G + a_{1}b_{1}$. Lemma \ref{lem t1p}(1) gives that $|D_{a_{1}b_{1}} \cap S| = 1$ and either $b_{1} \in D_{a_{1}b_{1}}$ or $a_{1} \in D_{a_{1}b_{1}}$. In the first case, $\{u, v\} \cap D_{a_{1}b_{1}} = \emptyset$ by Lemma \ref{lem 1}(3). To dominate $a_{2}$, $D_{a_{1}b_{1}} \cap (S \cap I) = \emptyset$. Hence, $D_{a_{1}b_{1}} \cap S = \emptyset$, a contradiction. Therefore, $a_{1} \in D_{a_{1}b_{1}}$. To dominate $H_{1} - b_{1}$ and by the connectedness of $(G+ a_{1}b_{1})[D_{a_{1}b_{1}}]$, $(D_{a_{1}b_{1}} - \{a_{1}\}) \subseteq \{u, v\}$. Lemma \ref{lem 1}(3) implies that $v \in D_{a_{1}b_{1}}$. Thus $v \succ H_{1} - b_{1}$. Let $b_{2} \in H_{1} - \{b_{1}\}$. Therefore $b_{2} \succ \{u, v\}$. Consider $G + a_{1}b_{2}$. Lemma \ref{lem t1p}(1) implies that we have $|D_{a_{1}b_{1}} \cap S| = 1$ and either $a_{1} \in D_{a_{1}b_{2}}$ or $b_{2} \in D_{a_{1}b_{2}}$. In the first case, $\{u, v\} \cap D_{a_{1}b_{2}} = \emptyset$ by Lemma \ref{lem 1}(3). By the connectedness of $(G + a_{1}b_{2})[D_{a_{1}b_{2}}]$, $(S \cap I) \cap D_{a_{1}b_{2}} = \emptyset$. Thus $D_{a_{1}b_{2}} \cap S = \emptyset$, a contradiction. Therefore, $b_{2} \in D_{a_{1}b_{2}}$. To dominate $a_{2}$, $(S \cap I) \cap D_{a_{1}b_{2}} = \emptyset$. Lemma \ref{lem 1}(3) yields that $D_{a_{1}b_{2}} \cap \{u, v\} = \emptyset$. Therefore $D_{a_{1}b_{2}} \cap S = \emptyset$, a contradiction and Case 2 cannot occur. So $p \geq 3$ and this completes the proof.
\qed

\indent In view of Lemma \ref{lem Co3p}, $p \geq 3$. Lemma \ref{lem 2} then implies that the vertices in $\cup^{2}_{i = 1}I_{i}$ can be ordered as $x_{1}, x_{2}, ..., x_{p}$ and there exists a path $y_{1}, y_{2}, ..., y_{p - 1}$ with $\{x_{i}, y_{i}\} \succ_{c} G + x_{i}x_{i + 1}$ for $i = 1, 2, ..., p - 1$.
\vskip 5 pt

\begin{lem}\label{lem Co4p}
$y_{i} \succ S \cap I$ and $y_{i} \in S - I$ for $i  = 1, 2, ..., p - 1$.
\end{lem}
\proof
Since $\{x_{i}, y_{i}\} \succ_{c} G + x_{i}x_{i + 1}$ for $i = 1, 2, ..., p - 1$ and $x_{i} \in I$, $y_{i} \succ S \cap I$. By the connectedness of $(G + x_{i}x_{i + 1})[D_{x_{i}x_{i + 1}}]$ and Lemma \ref{lem t1p}(1), $y_{i} \in S - I$ and this completes the proof.
\qed

\indent Lemma \ref{lem Co4p} implies that $\{y_{1}, y_{2}, ..., y_{p - 1}\} \subseteq S - I$. By Equation \ref{eq 2}, $|(S - I) - \{y_{1}, y_{2},$ $...,y_{p - 1}\}| = 1$. Let $\{y_{p}\} =(S - I) - \{y_{1}, y_{2}, ..., y_{p - 1}\}$.
\vskip 5 pt

\begin{lem}\label{lem Co5p}
For $i, j \in \{2, 3, ..., p\}$, if $y_{p}x_{i}, y_{p}x_{j} \in E(G)$, then $y_{i - 1}y_{j - 1} \in E(G)$.
\end{lem}
\proof
Consider $G + x_{i}x_{j}$. Lemma \ref{lem t1p}(1) yields that $|D_{x_{i}x_{j}} \cap \{x_{i}, x_{j}\}| = 1$ and $|D_{x_{i}x_{j}} \cap S| = 1$. Without loss of generality, let $x_{i} \in D_{x_{i}x_{j}}$ and $\{a\} = D_{x_{i}x_{j}} \cap S$. By the connectedness of $(G + x_{i}x_{j})[D_{x_{i}x_{j}}]$, $a \in S - I$. Since $x_{j} \succ (S - I) - \{y_{j - 1}\}$, it follows by Lemma \ref{lem 1}(3) that $a = y_{j - 1}$. Since $y_{i - 1}x_{i} \notin E(G)$, $y_{j - 1}y_{i - 1} \in E(G)$ and this completes the proof.
\qed
\vskip 5 pt

\begin{lem}\label{lem Co6p}
$\alpha_{1}, \alpha_{2} > 0$.
\end{lem}
\proof
By the assumption that $\alpha_{1} \leq \alpha_{2}$, we suppose that $\alpha_{1} = 0$. Clearly $\{x_{1}, x_{2}, ..., x_{p}\} \subseteq H_{2}$ and $y_{i} \succ H_{1}$ for $i = 1, 2, ..., p - 1$. Since $S$ is a minimum cut set, $N_{H_{1}}(y_{p}) \neq \emptyset$. Let $b \in N_{H_{1}}(y_{p})$. Therefore $b \succ S - I$. Consider $G + x_{1}b$. Lemma \ref{lem t1p}(1) implies that $|D_{x_{1}b} \cap S| = 1$ and either $b \in D_{x_{1}b}$ or $x_{1} \in D_{x_{1}b}$. Suppose that $b \in D_{x_{1}b}$. To dominate $x_{2}$, $D_{x_{1}b} \cap (S - I) \neq \emptyset$. Lemmas \ref{lem 2} and \ref{lem 1}(3) then imply that $D_{x_{1}b} \cap (S - I) = \{y_{p}\}$. So $y_{p} \succ \{x_{2}, x_{3}, ..., x_{p}\}$. Lemma \ref{lem Co5p} gives, further, that $G[y_{1}, y_{2}, ..., y_{p - 1}]$ is a clique. Lemma \ref{lem Co4p} then implies that $y_{i} \succ S \cap I$ for $i = 1, 2, ..., p - 1$. By Lemma \ref{lem Co1p}(2), $y_{i}y_{p} \notin E(G)$ for $i = 1, 2, ..., p - 1$. Therefore $y_{1}y_{p} \notin E(G)$. Because $\{x_{1}, y_{1}\} \succ_{c}G + x_{1}x_{2}$, $x_{1}y_{p} \in E(G)$, contradicting Lemma \ref{lem 1}(3). Therefore, $x_{1} \in D_{x_{1}b}$. By the connectedness of $(G + x_{1}b)[D_{x_{1}b}]$, $D_{x_{1}b} \cap (S \cap I) = \emptyset$. Lemma \ref{lem 1}(3) implies that $D_{x_{1}b} \cap (S - I) = \emptyset$. Thus $D_{x_{1}b} \cap S = \emptyset$ contradicting Lemma \ref{lem t1p}(1) and this completes the proof.
\qed
\vskip 5 pt

\begin{thm}\label{thm Co7p}
Let $G$ be a maximal $3$-$\gamma_{c}$-vertex critical graph. If $\alpha = \kappa$, then $G - S$ contains at least one singleton component for all minimum cut set $S$.
\end{thm}
\proof
Let a graph $G$ be a maximal $3$-$\gamma_{c}$-vertex critical graph with $\alpha = \kappa$. By Equation \ref{eq 2}, $|S - I| = \alpha_{1} + \alpha_{2}$. Suppose that there is no singleton component of $G - S$, in particular, $|H_{i}| > 1$ for $i = 1, 2$. Let $\alpha_{1} + \alpha_{2} = p$. Lemma \ref{lem Co3p} implies that $p \geq 3$. Lemma \ref{lem Co6p} gives also that $0 < \alpha_{1} \leq \alpha_{2}$. We also define $x_{1}, x_{2}, ..., x_{p}$ a path $y_{1}, y_{2}, ..., y_{p - 1}$ and a vertex $y_{p}$ by the same as in the previous lemmas.
\vskip 5 pt

\indent Suppose there exist $x_{i}, x_{j}$ for $i, j \in \{2, 3, ..., p\}$ such that $y_{p} \in D_{x_{i}x_{j}}$. Lemma \ref{lem 1}(1) and (2) then implies that $D_{x_{i}x_{j}} = \{x_{i}, y_{p}\}$ or $D_{x_{i}x_{j}} = \{x_{j}, y_{p}\}$. Without loss of generality, let $D_{x_{i}x_{j}} = \{x_{j}, y_{p}\}$. Thus $y_{p} \succ \{x_{1}, x_{2}, ..., x_{p}\} - \{x_{i}\}$. Since $\{x_{i}, y_{i}\} \succ_{c} G + x_{i}x_{i + 1}$, $y_{i}y_{p} \in E(G)$. Lemma \ref{lem Co5p} yields that $G[\{y_{1}, y_{2}, ..., y_{p - 1}\} - \{y_{i - 1}\}]$ is a clique. Since $y_{i}y_{i - 1} \in E(G)$, $y_{i} \succ S - I$. Lemma \ref{lem Co4p} implies that $y_{i} \succ S \cap I$. Therefore $y_{i} \succ S$, contradicting Lemma \ref{lem Co1p}(2). Hence, $y_{p} \notin D_{x_{i}x_{j}}$ for any $i, j \in \{2, 3, ..., p\}$. By the same arguments as in the proof of Lemma \ref{lem Co5p}, $G[\{y_{1}, y_{2}, ..., y_{p - 1}\}]$ is a clique. As $y_{i} \succ S \cap I$, by Lemma \ref{lem Co1p}(2), we must have $y_{i}y_{p} \notin E(G)$ for $i \in \{1, 2, ..., p - 1\}$. Since $\{x_{i}, y_{i}\} \succ_{c} G + x_{i}x_{i + 1}$ for $i \in \{1, 2, ..., p - 1\}$, $x_{i}y_{p} \in E(G)$. So $x_{1} \succ S - I$. By Lemma \ref{lem Co1p}(2), $S \cap I \neq \emptyset$ as otherwise $x_{1} \succ S$. Let $x_{1} \in H_{i}$ for some $i \in \{1, 2\}$. Consider $G - x_{1}$. Since $|H_{j}| > 1$ for $j = 1, 2$, neither $D_{x_{1}} \subseteq H_{1}$ nor $D_{x_{1}} \subseteq H_{2}$. Lemma \ref{lem Co1p}(1) gives, further, that $D_{x_{1}} \cap S \neq \emptyset$. Lemma \ref{lem 6}(2) implies that $D_{x_{1}} \cap (S - I) = \emptyset$. Thus $D_{x_{1}} \cap (S \cap I) \neq \emptyset$. Let $u_{1} \in D_{x_{1}} \cap (S \cap I)$. Lemma \ref{lem 6}(1) implies also that $|D_{x_{1}} - \{u_{1}\}| = 1$. Let $\{w\} = D_{x_{1}} - \{u_{1}\}$. If $w \in H_{i}$, then $u_{1} \succ H_{3 - i}$. Since $u_{1} \in I$, $\alpha_{3 - i} = 0$ contradicting Lemma \ref{lem Co6p}. So $w \in H_{3 - i}$ and $u_{1} \succ H_{i} - x_{1}$. Since $u_{1} \in I$, $I_{i} = \{x_{1}\}$. It follows that $\{x_{2}, x_{3}, ..., x_{p}\} \subseteq H_{3 - i}$.
\vskip 5 pt

\noindent \textbf{Claim 1 :} There is no $u \in S \cap I$ such that $u \succ S - I$.\\
\indent Suppose $u \succ S - I$. Consider $G - u$. Lemma \ref{lem Co1p}(1) implies that $D_{u} \cap S \neq \emptyset$. Lemma \ref{lem 6}(2) implies also that $D_{u} \cap (S - I) = \emptyset$. Thus there exists $u' \in D_{u} \cap (S \cap I)$. Lemma \ref{lem 6}(1) gives that $|D_{u} - \{u'\}| = 1$. Let $\{z\} = D_{u} - \{u'\}$. To dominate $x_{1}$, $z \in H_{i}$. So $D_{u}$ does not dominate $I_{3 - i}$, a contradiction and thus establishing Claim 1.
\vskip 5 pt

\indent Claim 1 together with Lemma \ref{lem Co4p} yield that $y_{p}$ is not adjacent to any vertex in $S \cap I$. We now have $y_{p}$ an isolated vertex in $S$.
\vskip 5 pt

\noindent \textbf{Claim 2 :} $y_{1} \succ H_{i}$.\\
\indent Suppose $y_{1}$ is not adjacent to $b_{1} \in H_{i}$. Consider $G + b_{1}x_{2}$. We see that $b_{1}y_{1}, x_{2}y_{1} \notin E(G)$. Lemma \ref{lem t1p}(1) gives that $|D_{b_{1}x_{2}} \cap S| = 1$and either $b_{1} \in D_{b_{1}x_{2}}$ or $x_{2} \in D_{b_{1}x_{2}}$. If $b_{1} \in D_{b_{1}x_{2}}$, then $(S - \{y_{1}, y_{p}\}) \cap D_{b_{1}x_{2}} = \emptyset$ to dominate $I_{3 - i}$. Since $y_{p}x_{2} \in E(G)$, by Lemma\ref{lem 1}(3), $y_{p} \notin D_{b_{1}x_{2}}$. By the connectedness of $(G + b_{1}x_{2})[D_{b_{1}x_{2}}]$, $y_{1} \notin D_{b_{1}x_{2}}$. Therefore $D_{b_{1}x_{2}} \cap S = \emptyset$, a contradiction. Hence, $x_{2} \in D_{b_{1}x_{2}}$. To dominate $I_{3 - i} \cup (S \cap I)$, $D_{b_{1}x_{2}} \cap \{y_{2}, y_{3}, ..., y_{p}\} = \emptyset$. By the connectedness of $(G + b_{1}x_{2})[D_{b_{1}x_{2}}]$, $((S \cap I) \cup \{y_{1}\}) \cap D_{b_{1}x_{2}} = \emptyset$. Therefore, $D_{b_{1}x_{2}} \cap S = \emptyset$, a contradiction and we settle Claim 2.
\vskip 5 pt

\indent Let $b_{1} \in H_{i} - \{x_{1}\}$. Recall that $u_{1} \succ H_{i} - x_{1}$. Clearly $b_{1}u_{1} \in E(G)$. By Claim 2 and Lemma \ref{lem 2}, $b_{1} \succ \{y_{1}, y_{2}, ..., y_{p - 1}\} \cup \{u_{1}\}$. Consider $G - b_{1}$. Lemma \ref{lem Co1p}(1) implies that $D_{b_{1}} \cap S \neq \emptyset$. Lemma \ref{lem 6}(2) gives that $D_{b_{1}} \cap (\{y_{1}, y_{2}, ..., y_{p - 1}\} \cup \{u_{1}\}) = \emptyset$. If there is $u_{2} \in D_{b_{1}} \cap ((S \cap I) - \{u_{1}\})$, then, by Lemma \ref{lem 6}(1), let $\{y'\} = D_{b_{1}} - \{u_{2}\}$. To dominate $x_{1}$, $y' \in H_{i}$. Thus $D_{b_{1}}$ does not dominate $x_{2}$, a contradiction. Therefore, $\{y_{p}\} = D_{b_{1}} \cap S$. Since $y_{p}$ is an isolated vertex in $S$, by Lemma \ref{lem t1p}(2), at least one of $C_{i}$ is a singleton component contradicting the assumption. We then finish the proof.
\qed

\indent By Theorem \ref{thm Co7p}, we obtain the following corollary.
\begin{cor}\label{cor p}
Let $G$ be a maximal $3$-$\gamma_{c}$-vertex critical graph. If $\alpha = \kappa$, then $\kappa = \delta$.
\end{cor}
\proof
By Theorem \ref{thm Co7p}, we have at least one singleton component $C_{i}$ of $G - S$. Let $\{c\} = V(C_{i})$. Thus $N_{G}(c) \subseteq S$. Hence, $\delta \leq deg_{G}(c) \leq |S| = \kappa \leq \delta$.
\qed

\indent We next construct the class $\mathcal{G}_{4}(s)$ of maximal $3$-$\gamma_{c}$-vertex critical graphs which $\alpha < \kappa$ and $\kappa < \delta$ in order to show that the condition $\alpha = \kappa$ is necessary to prove Corollary \ref{cor p}. Let $R, T, W$ and $Z$ be non-empty disjoint vertex sets where $R = \{r_{1}, r_{2}, ..., r_{s}\}, T = \{t_{1}, t_{2}, ..., t_{s}\}, W = \{w_{1}, w_{2}, ..., w_{s}\}, Z = \{z_{1}, z_{2}, ..., z_{s}\}$ and $s \geq 3$.  A graph $G$ in the class $\mathcal{G}_{4}(s)$ can be constructed from $R, T, W$ and $Z$ by adding edges according the join operations :
\begin{itemize}
  \item for $1 \leq i \leq s$, $r_{i} \vee R \cup T \cup W - \{r_{i}, t_{i}\}$,
  \item $t_{i} \vee R \cup W \cup Z - \{w_{i}, r_{i}\}$,
  \item $w_{i} \vee R \cup T \cup Z - \{t_{i}, z_{i}\}$,
  \item $z_{i} \vee Z \cup T \cup W = \{z_{i}, w_{i}\}$ and
  \item adding edges so the the vertices in $R$ and $Z$ form cliques.
\end{itemize}
\noindent Observe that, for $1 \leq i \leq s$, $N_{G}(r_{i}) = R \cup T \cup W - \{r_{i}, t_{i}\}$, $N_{G}(t_{i}) = R \cup W \cup Z - \{w_{i}, r_{i}\}$, $N_{G}(w_{i}) = R \cup T \cup Z - \{t_{i}, z_{i}\}$ and $N_{G}(z_{i}) = Z \cup T \cup W = \{z_{i}, w_{i}\}$. Observe also that $T$ and $W$ are independent sets. A graph $G$ is illustrated in Figure 4 where double lines between two sets means joining each vertex in one set to every vertex in the other set.
\vskip 20 pt

\setlength{\unitlength}{0.7cm}
\begin{picture}(0, 9)

\put(3.2, 6.5){\circle*{0.2}}
\put(4.2, 6.5){\circle*{0.2}}
\put(6.8, 6.5){\circle*{0.2}}

\put(3.2, 2.05){\circle*{0.2}}
\put(4.2, 1.45){\circle*{0.2}}
\put(6.8, 0.5){\circle*{0.2}}

\put(10.2, 0.5){\circle*{0.2}}
\put(12.8, 1.45){\circle*{0.2}}
\put(13.8, 2.05){\circle*{0.2}}

\put(10.2, 6.5){\circle*{0.2}}
\put(12.8, 6.5){\circle*{0.2}}
\put(13.8, 6.5){\circle*{0.2}}

\multiput(3.2, 2.05)(0,0.65){7}{\line(0,1){0.4}}
\multiput(4.2, 1.45)(0,0.65){8}{\line(0,1){0.4}}
\multiput(6.8, 0.5)(0,0.6){10}{\line(0,1){0.4}}

\multiput(10.2, 0.5)(0,0.6){10}{\line(0,1){0.4}}
\multiput(12.8, 1.45)(0,0.65){8}{\line(0,1){0.4}}
\multiput(13.8, 2.05)(0,0.65){7}{\line(0,1){0.4}}

\multiput(3.2, 2.05)(0.6, 0){18}{\line(1,0){0.4}}
\multiput(4.4, 1.45)(0.6, 0){14}{\line(1,0){0.4}}
\multiput(6.8, 0.5)(0.6, 0){6}{\line(1,0){0.4}}

\put(5, 6.5){\oval(5, 1.7)}
\put(12, 6.5){\oval(5, 1.7)}

\put(2.5,0){\dashbox{0.2}(5,2.5)[s]}
\put(9.5,0){\dashbox{0.2}(5,2.5)[s]}

\put(5, 2.5){\line(0, 1){3.15}}
\put(5.1, 2.5){\line(0, 1){3.15}}

\put(12, 2.5){\line(0, 1){3.15}}
\put(11.9, 2.5){\line(0, 1){3.15}}

\put(7.5, 0.8){\line(1, 0){2}}
\put(7.5, 0.9){\line(1, 0){2}}

\put(7.5, 2.5){\line(1, 2){2}}
\put(7.5, 2.35){\line(1, 2){2}}

\put(9.5, 2.5){\line(-1, 2){2}}
\put(9.5, 2.35){\line(-1, 2){2}}
\put(11.2, 1){$...$}
\put(11.2, 6.5){$...$}
\put(5.3, 1){$...$}
\put(5.3, 6.5){$...$}

\put(3.2, 6.7){$r_{1}$}
\put(4.2, 6.7){$r_{2}$}
\put(6.8, 6.7){$r_{s}$}

\put(3, 1.55){$t_{1}$}
\put(4, 0.95){$t_{2}$}
\put(6.6, 0){$t_{s}$}

\put(14, 1.55){$w_{1}$}
\put(13, 0.95){$w_{2}$}
\put(10.4, 0){$w_{s}$}

\put(13.8, 6.7){$z_{1}$}
\put(12.8, 6.7){$z_{2}$}
\put(10.2, 6.7){$z_{s}$}

\put(5, 7.5){$R$}
\put(5, -0.5){$T$}
\put(12, -0.5){$W$}
\put(12, 7.5){$Z$}
\put(5, -1.9){\footnotesize\textbf{Figure 4} : A graph $G$ in the class $\mathcal{G}_{4}(s)$}
\end{picture}
\vskip 70 pt

\begin{lem}\label{lem g4}
If $G \in \mathcal{G}_{4}(s)$, then $G$ is a maximal $3$-$\gamma_{c}$-vertex critical graph.
\end{lem}
\proof We see that $\{r_{1}, t_{2}, w_{2}\} \succ_{c} G$. Thus $\gamma_{c}(G) \leq 3$. Suppose there exist $u, v \in V(G)$ such that $\{u, v\} \succ_{c} G$. Let $i \in \{1, .., s\}$. Suppose that $u = r_{i}$. To dominate $Z$, $v \notin R$. If $v \in T$, then, by connected, $v \neq t_{i}$. Thus $\{u, v\}$ does not dominate $t_{i}$. To dominate $Z$, $v \notin W$. Clearly $v \in Z$ but $G[\{u, v\}]$ is not connected, a contradiction. Hence, $\{u, v\} \cap R = \emptyset$. By symmetric, $\{u, v\} \cap Z = \emptyset$. Therefore $\{u, v\}\subseteq T \cup W$. Without loss of generality, $u = t_{i}$. By connected, $v \in W - \{w_{i}\}$. Thus $\{u, v\}$ does not dominate $w_{i}$. Therefore $\gamma_{c}(G) = 3$.
\vskip 5 pt

\indent To establish the criticality, let $u, v$ be a pair of non-adjacent vertices of $G$. For $1 \leq i \leq s$, if $\{u, v\} = \{r_{i}, t_{i}\}$, then $D_{uv} = \{r_{i}, t_{i}\}$. If $\{u, v\} = \{t_{i}, w_{i}\}$, then $D_{uv} = \{t_{i}, w_{i}\}$. If $\{u, v\} = \{w_{i}, z_{i}\}$, then $D_{uv} = \{w_{i}, z_{i}\}$. For $1 \leq i \neq j \leq s$, if $\{u, v\} = \{t_{i}, t_{j}\}$, then $D_{uv} = \{t_{i}, r_{j}\}$. If $\{u, v\} = \{w_{i}, w_{j}\}$, then $D_{uv} = \{w_{i}, z_{j}\}$. If $\{u, v\} = \{r_{i}, z_{l}\}$ where $l \in \{1, 2, ..., s\}$, then $D_{uv} = \{r_{i}, z_{l}\}$. Thus $G$ is a $3$-$\gamma_{c}$-edge critical graph. Let $v \in V(G)$. For $1 \leq i \neq j \leq s$, if $u = r_{i}$, then $D_{v} = \{t_{i}, z_{j}\}$. If $v = t_{i}$, then $D_{v} = \{t_{j}, r_{i}\}$. If $v = w_{i}$, then $D_{v} = \{z_{i}, w_{j}\}$. Finally, if $v = z_{i}$, then $D_{v} = \{w_{i}, r_{j}\}$. Therefore $G$ is a maximal $3$-$\gamma_{c}$-vertex critical graph and this completes the proof.
\qed
\vskip 5 pt

\indent We see that $G$ has $T$ a maximum independent set and $T \cup W$ a mimimum cut set. Thus $\alpha = s < 2s = \kappa$. We, further, have that $G$ is a regular graph with $deg_{G}(v) = 3s - 2$ for all $v \in V(G)$. Since $s \geq 3$, $\delta = 3s - 2 > 2s = \kappa$. Hence, the condition $\alpha = \kappa$ is necessary to prove Corollary \ref{cor p}.
\vskip 5 pt

\indent We conclude this paper by establishing hamiltonian property of maximal $3$-$\gamma_{c}$-vertex critical graphs. We prove the following result by using Theorem \ref{thm CE}.
\vskip 5 pt

\begin{cor}\label{cor Co8p}
For any $3$-connected maximal $3$-$\gamma_{c}$-vertex critical graph $G$. If $\kappa < \delta$, then $G$ is hamiltonian connected.
\end{cor}
\proof
Suppose that $\kappa < \delta$. Theorem \ref{thm Co2p} and Corollary \ref{cor p} then imply that $\alpha < \kappa$. By Theorem \ref{thm CE}, $G$ is hamiltonian connected.
\qed
\vskip 5 pt

\indent However, we are sure that every $3$-connected maximal $3$-$\gamma_{c}$-vertex critical graph is hamiltonian connected. To prove this, by Theorem \ref{thm CE} and Corollary \ref{cor Co8p},we need only prove the following conjecture.
\vskip 5 pt

\begin{conj}
Let $G$ be a $3$-connected maximal $3$-$\gamma_{c}$-vertex critical graph with $\alpha = \kappa = \delta$. Then $G$ is hamiltonian connected.
\end{conj}



\begin{thebibliography}{99}\label{bib}
\bibitem{A} N. Ananchuen, On Domination Critical Graphs with Cut Vertices having Connected Domination Number 3, International Mathematical Forum 2(2007) 3041 - 3052.


\bibitem{AAP} W. Ananchuen, N. Ananchuen and M.D. Plummer, Vertex Criticality for Connected Domination, Utilitas Mathematica 86(2011) 45 - 64.


\bibitem{BM} J. A. Bondy and U. S. R. Murty, Graphs Theory with Applications, The Macmillan Press, London, 1976.





\bibitem{CSM} X.G. Chen, L.Sun and D.X. Ma, Connected Domination Critical Graphs, Applied Mathematics Letters 17(2004) 503-507.

\bibitem{CE} V. Chv\'{a}tal and P. Erd\"{o}s, A Note on Hamiltonian Cirtcuits, Discrete Mathematics 2(1972) 111 - 113.




\bibitem{KC} P. Kaemawichanurat and L. Caccetta, Independence and Connectivity of Connected Domination Critical Graphs, manuscript.

\bibitem{KCA} P. Kaemawichanurat, L. Caccettaa  and W. Ananchuen, Hamiltonicity of Connected Domination Critical Graphs. submitted for publication.

\bibitem{KCA1} P. Kaemawichanurat, Connected Domination Critical Graphs. Ph.D. Thesis. PhD supervisor: Louis Caccetta. Curtin University (2016).

\bibitem{KCA2} P. Kaemawichanurat, On The Independence Number of $3$-$\gamma_{c}$-Edge Critical Graphs.

\bibitem{S} J. Simmons, Closure Operations and Hamitonian Properties of Independent and Total Domination Critical Graphs. Ph.D. Thesis. PhD advisor: Gary MacGillvray. University of Victoria(2005).

\bibitem{ZT} L. Zhang and F. Tian, Independence and Connectivity in $3$-Domination Critical Graphs, Discrete Mathematics 259(2002) 227 - 236.

\end{thebibliography}
\end{document}